\documentclass[12pt,draft]{amsart}
\title{Maximality of the microstates free entropy for $R$-diagonal
  elements}
\author{Alexandru Nica}
\thanks{Research supported by a grant from the Natural Sciences 
and Engineering Research Council, Canada. } 
\address{
Department of Pure Mathematics,
University of Waterloo,
Waterloo, Ontario N2L 3G1, Canada}
\email{anica@math.uwaterloo.ca} 
\author{Dimitri Shlyakhtenko}
\thanks{Research supported in part by a National Science Foundation
  postdoctoral fellowship DMS-9804625}
\address{Department of Mathematics, UCLA, Los Angeles, CA 90095}
\email{shlyakht@member.ams.org}
\author{Roland Speicher}
\thanks{Supported by a Heisenberg Fellowship of the DFG, Germany. } 
\address{
Institut f\"ur Angewandte Mathematik,
Universit\"at Heidelberg,
D-69120 Heidelberg, Germany}
\email{roland.speicher@urz.uni-heidelberg.de } 
\date\today
\newcommand{\sa}{{\operatorname{sa}}}
\newcommand{\chisa}{{\chi^\sa}}

\newcommand{\GammasaR}{{\Gamma^\sa_R}}
\newcommand{\Tr}{\operatorname{Tr}}
\renewcommand{\Re}{\operatorname{Re}}
\renewcommand{\Im}{\operatorname{Im}}
\newcommand{\FF}{Q}
\newcommand{\FFF}{T}
\newtheorem{theorem}[subsection]{Theorem}
\newtheorem{lemma}[subsection]{Lemma}
\newtheorem{proposition}[subsection]{Proposition}%
\newtheorem{definition}[subsection]{Definition}%
\newtheorem{theorem*}{Theorem}%
\theoremstyle{remark}
\newtheorem{remark}[subsection]{Remark}%
\newtheorem{notation}[subsection]{Notation}%
\begin{document}%
\begin{abstract}
  An non-commutative non-self adjoint random variable $z$ is called
  $R$-diagonal, if its $*$-distribution is invariant under
  multiplication by free unitaries: if a unitary $w$ is $*$-free from
  $z$, then the $*$-distribution of $z$ is the same as that of $wz$.
  Using Voiculescu's microstates definition of free entropy, we show
  that the $R$-diagonal elements are characterized as having the
  largest free entropy among all variables $y$ with a fixed
  distribution of $y^*y$.  More generally, let $Z$ be a $d\times d$
  matrix whose entries are non-commutative random variables $X_{ij}$,
  $1\leq i,j\leq d$.  Then the free entropy of the family
  $\{X_{ij}\}_{ij}$ of the entries of $Z$ is maximal among all $Z$
  with a fixed distribution of $Z^*Z$, if and only if $Z$ is
  $R$-diagonal and is $*$-free from the algebra of scalar $d\times d$
  matrices.  The results of this paper are analogous to the results of
  our paper \cite{nss:entropy}, where we considered the same problems
  in the framework of the non-microstates definition of entropy.
\end{abstract}
\maketitle
\section{Introduction.}
Let $(M,\tau)$ be a tracial non-commutative $W^*$-probability space.
A (non-self-adjoint) element $z\in M$ is called $R$-diagonal if its
$*$-distribution is invariant under multiplication by free unitaries;
i.e., if $u$ is a unitary, $*$-free from $z$, the $*$-distributions of
$uz$ and $z$ coincide.  The concept of $R$-diagonality was introduced
in \cite{speicher:r-diag}, where it was shown to be equivalent to
several conditions; we mention that if $z^*z$ has a (possibly
unbounded) inverse (in particular, if the distribution of $z^*z$ is
non-atomic), then $z$ is $R$-diagonal if and only if in its polar
decomposition $z=u(z^*z)^{1/2}$, $u$ is $*$-free from $(z^*z)^{1/2}$
and satisfies $\tau(u^k)=0$ for $k\in\mathbb{Z}\setminus \{0\}$.

In our recent paper \cite{nss:entropy} $R$-diagonal elements appeared
in connection with certain maximization problems in free entropy.
Free entropy was introduced by Voiculescu in \cite{dvv:entropy2};
later, a different definition was given by him in \cite{dvv:entropy5}.
The first definition involves approximating the given $n$-tuple of
variables using finite-dimensional matrices (so-called microstates);
the normalized limit of the logarithms of volumes of all such possible
microstates is then the free entropy.  On the other hand, Voiculescu's
definition in \cite{dvv:entropy5} does not involve microstates, but
uses free Fisher information measure and non-commutative Hilbert
transform.  At present it is not known whether the two definitions of
free entropy always give the same quantity.  Our approach in
\cite{nss:entropy} used the second definition of Voiculescu.

In this paper we prove two theorems for the microstates free entropy,
which are analogous to our results in \cite{nss:entropy} for the
second (non-microstates) definition of entropy.  One of our results
can be interpreted as saying that $R$-diagonal elements $z$ are
characterized by the statement that the free entropy $\chi(z)$ is
maximal among all possible $\chi(y)$, so that the distributions of
$y^*y$ and $z^*z$ are the same.

When this paper was almost finished we received a preprint of Hiai and
Petz \cite{hiai-petz}, where the same kind of problems were
considered.

If $Y_1,\dots,Y_n\in M$ (not necessarily self-adjoint), we denote by
$\chi(Y_1,\dots,Y_n)$ the free entropy of $Y_1,\dots,Y_n$ as defined
by Voiculescu in \cite{dvv:improvedrandom}.  We denote by
$\chisa(X_1,\dots,X_n)$ for $X_i\in M$ self-adjoint the free entropy
of a self-adjoint $n$-tuple as defined in \cite{dvv:entropy2}; we give
a brief review of these quantities below in \S\ref{sect:defs}. A
unitary $u$ in a non-commutative probability space $(M,\tau)$ is
called a Haar unitary if $\tau(u^k)=0$ for all
$k\in\mathbb{Z}\setminus\{0\}$.

\begin{theorem*} 
  Let $y\in M$, and let $u\in M$ be a Haar unitary which is $*$-free
  from $b=(y^*y)^{1/2}$.  Let $x$ be an element such that
  $\tau(x^{2k})=\tau(b^{2k})$ and $\tau(x^{2k+1})=0$, for all
  $k\in\mathbb N$ (i.e., $x$ is symmetric). Then
  \begin{enumerate} 
    \renewcommand{\theenumi}{\arabic{section}.\arabic{subsection}(\alph{enumi})}
    \renewcommand{\labelenumi}{(\alph{enumi})}
  \item $\chi(y) \leq \chi(ub)$.
  \item $\chi(ub) = \chisa(b^2/2)+3/4+1/2 \log 2\pi = 2
    \chisa(2^{-\frac{1}{2}} x)$
  \item If $\chi(y)=\chi(ub)>-\infty$, then $y$ is $R$-diagonal, i.e.,
    in the polar decomposition $y=vb$ we have: $v$ is a Haar unitary
    and is $*$-free from $b$.
  \end{enumerate}
\end{theorem*}

Let $\omega\in\beta \mathbb{N}\setminus\mathbb{N}$ be a free
ultrafilter; i.e., a homomorphism from the algebra $C(\mathbb{N})$ of
all bounded (continuous) functions on $\mathbb N$ to $\mathbb C$,
which is not given by the evaluation at a point in $\mathbb N$.  For
$d\in\mathbb N$ we write $d\omega$ for the free ultrafilter
corresponding to the functional $f\mapsto \lim_{n\to\omega} f(dn)$.
Given $\omega$, one can construct (see \cite{dvv:improvedrandom} and
see also a brief review below) free entropy quantities
${\chi^\sa}^\omega$ and $\chi^\omega$, which have properties similar
to those of $\chi^\sa$ and $\chi$; it is in fact not known whether
these quantities are different.  It is known that in the one-variable
case, $\chi^\sa(X)={\chi^\sa}^\omega(X)$.

\begin{theorem*}
  Let $X_{ij}$, $1\leq i,j\leq d$ be a family of non-commutative
  random variables in a tracial non-commutative probability space
  $(M,\hat\tau)$.  Let $Z\in M\otimes M_d$ be given by
  $$Z=\sum_{i,j=1}^d X_{ij}\otimes e_{ij},$$ where $e_{ij}$ are matrix
  units in the algebra of $d\times d$ matrices. We denote by $\tau$
  the normalized trace on $M\otimes M_d$.  Let $\omega$ be a free
  ultrafilter.  Let $X$ be a self-adjoint variable with
  $\tau(X^{2n+1})=0$ for all $n\in\mathbb N$, and such that $\tau(X^{2n})
  =\tau((Z^*Z)^n)$, $\forall n\in\mathbb N$.  Then we have
  \begin{enumerate}
    \renewcommand{\theenumi}%
    {\arabic{section}.\arabic{subsection}(\alph{enumi})}
    \renewcommand{\labelenumi}{(\alph{enumi})}
  \item $\chi^\omega(\{X_{ij}\}_{1\leq i,j\leq d})\leq d^2
    \chi^{d\omega}(Z) + {d^2} \log d\leq 2d^2
    \chi^\sa(2^{\frac{-1}{2}}X) + {d^2} \log d$;
  \item Equality holds in (a) if $Z$ is $R$-diagonal and $*$-free from
    the algebra $1\otimes M_d$.
  \item If equality holds in (a) and $\chi^\sa(2^{\frac{-1}{2}}X)\neq
    -\infty$, then $Z$ is $R$-diagonal and is $*$-free from the algebra
    $1\otimes M_d$.
  \end{enumerate}
\end{theorem*}

The proof of the first theorem is quite different in nature than our
proof in \cite{nss:entropy} (the microstates-free proof relied on the
notion of free entropy with respect to a completely-positive map
introduced in \cite{shlyakht:cpentropy}). On the other hand, the proof
of the second theorem is analogous to the one we gave in
\cite{nss:entropy}, and relies on the microstates analog
\cite{shlyakht:relmicro} of the relative entropy \cite{dvv:entropy5}
that we used in the microstates-free approach.

\section{Maximality of microstates free entropy for $R$-diagonal
  pairs} Let $(M,\tau)$ be a tracial $W^*$-probability space, and
$b\in M$ be a fixed positive element. Let $u\in M$ be a Haar unitary
which is $*$-free from $b$.  Lastly, let $x\in M$ be such that for all
$k\in\mathbb N$, $\tau(x^{2k+1})=0$ and $\tau(x^{2k})=\tau(b^{2k})$.
The main result of the section is

\begin{theorem} \label{thrm:microstates}
  Let $u$, $b$ and $x$ be as above.  Assume that $y\in M$ satisfies
  $(y^*y)^{1/2}=b$.  Then
  \begin{enumerate} 
    \renewcommand{\theenumi}{\arabic{section}.\arabic{subsection}(\alph{enumi})}
    \renewcommand{\labelenumi}{(\alph{enumi})}
  \item $\chi(y) \leq \chi(ub)$. \label{thrm:micro:1}
  \item $\chi(ub) = \chisa(b^2/2)+3/4+1/2 \log 2\pi = 2
    \chisa(2^{-\frac{1}{2}} x)$ \label{thrm:micro:2}
  \item If $\chi(y)=\chi(ub)>-\infty$, then $y$ is $R$-diagonal, i.e.,
    in the polar decomposition $y=vb$, we have: $v$ is a Haar unitary
    and is $*$-free from $b$. \label{thrm:micro:3}
  \end{enumerate}
  The same conclusions hold for $\chi^\omega$ in place of $\chi$.
\end{theorem}

Before starting the proof of the theorem, we fix some notation and
definitions.

\begin{notation}  We use the following notation 
\begin{itemize}
\item $U(k)$ is the unitary group of $k\times k$ unitary matrices.
\item $M_k$ is the set of all $k\times k$ matrices; $M_k^\sa$ is the
  set of all self-adjoint matrices in $M_k$.
\item $M_k^+\subset M_k$ is the set of all positive $k\times k$ matrices.
\item $\mu_k$ is the normalized Haar measure on $U(k)$; thus
  $\mu_k(U(k))=1$.
\item $\lambda_k$ is the measure on $M_k$, coming from its Euclidean
  structure $\langle a,b\rangle =\Re\Tr(ab^*)$, where $\Tr$ is the
  usual matrix trace, $\Tr(I)=k$; $\lambda_k^\sa$ is the Lebesgue
  measure on $M_k^\sa$ coming from its Euclidean structure $\langle
  a,b\rangle =\Re\Tr(ab^*)$.
\item $\lambda_k^+$ is the measure on $M_k^+$ coming from its structure of
  a cone in the Euclidean space of $k\times k$ matrices.
\item $P:U(k)\times M_k^+\to M_k$ is given by $(v,p)\mapsto vp$
\item $\Omega_k$ is the canonical volume form on $M_k$ giving rise to
  Lebesgue measure. 
\item $\Omega^u_k \wedge \Omega^+_k$ is the canonical volume form on
  $U(k)\times M_k^+$, giving rise to the product measure $\mu_k\times
  \lambda^+_k$.
\item $\mathfrak{u}(k)$ is the Lie algebra of $U(k)$. 
\item $C_k$ is the volume of $U(k)$ with respect to the bi-invariant
  volume form arising from the Euclidean structure on ${\mathfrak
    u}(k)$ coming from the Killing form $\langle a,b\rangle =
  \Re\Tr(ab)$.
\end{itemize}
\end{notation}

\subsection{Definitions of free entropy.}  \label{sect:defs}
Let $X_1,\dots,X_n\in M$ be self-adjoint, and $Y_1,\dots,Y_n\in M$ be
not necessarily self-adjoint.  Let $\epsilon>0$, $R>0$ be real numbers
and $k>0$, $m>0$ be integers.  Then define the sets (cf.
\cite{dvv:entropy2,dvv:improvedrandom})
\begin{eqnarray*}
  \Gamma^\sa_R(X_1,\dots,X_n;m,k,\epsilon)&=&\{(x_1,\dots,x_n)\in
  (M_k^\sa)^n : \\ & & |\frac{1}{k}\Tr(x_{i_1}\dots x_{i_p})-
  \tau(X_{i_1}\dots X_{i_p})|<\epsilon\\ & & \text{for all } p\leq m,
  1\leq i_j\leq n, 1\leq j\leq p \};\\ 
  \Gamma_R(Y_1,\dots,Y_n;m,k,\epsilon)&=&\{(y_1,\dots,y_n)\in (M_k)^n:
  \\ & & |\frac{1}{k}\Tr(y^{g_1}_{i_1}\dots y^{g_p}_{i_p})-
  \tau(Y_{i_1}^{g_1} \dots Y_{i_p}^{g_p})|<\epsilon\\ & & \text{for all } p\leq m,
  1\leq i_j\leq n, g_j\in \{*,\cdot\}, 1\leq j\leq p \};
\end{eqnarray*}

Define next $$\chi^\sa(X_1,\dots,X_n;m,\epsilon)=\limsup_{k\to\infty}
\left[\frac{1}{k^2}\log\lambda_k
\Gamma_R^\sa(X_1,\dots,X_n;m,k,\epsilon) + \frac{n}{2} \log k\right]$$
and similarly
$$\chi(Y_1,\dots,Y_n;m,\epsilon)=\limsup_{k\to\infty}
\left[\frac{1}{k^2}\log\lambda_k \Gamma_R(Y_1,\dots,Y_n;m,k,\epsilon)
+ {n} \log k\right].$$ For $\omega$ a free ultrafilter on $\mathbb N$,
the quantities $\chi^\omega(Y_1,\dots,Y_n;m,\epsilon)$ and
${\chi^\sa}^\omega(X_1,\dots,X_n;m,\epsilon)$ are defined in exactly
the same way, except that $\limsup_{k\to\infty}$ is replaced by
$\lim_{k\to\omega}$.  Next, the free entropy is defined by
$$\chi^\sa(X_1,\dots,X_n)=\sup_R \inf_{m,\epsilon}
\chi^\sa(X_1,\dots,X_n;m,\epsilon);$$ the quantities
${\chi^\sa}^\omega$, $\chi$, $\chi^\omega$ are defined in exactly the
same way, using in the place of $\chi^\sa(\cdots;m,\epsilon)$ the
quantities ${\chi^\sa}^\omega(\cdots;m,\epsilon)$,
$\chi(\cdots;m,\epsilon)$, and $\chi^\omega(\cdots;m,\epsilon)$,
respectively.

\begin{definition}
  Let $(X_R(k,m,\epsilon),\mu^X_{R,k,m,\epsilon})$ and
  $(Y_R(k,m,\epsilon),\mu^Y_{R,k,m,\epsilon})$ be two sequences of
  measure spaces depending on $k,m\in\mathbb N$ and $R,\epsilon\in
  (0,+\infty)$.  We shall say that $X$ is asymptotically included in
  $Y$, if for all $m$, $\epsilon$, $R$, there is $k_0$, $m'\geq m$,
  $\epsilon'\leq \epsilon$, $R'>R$, such that for all $k>k_0$, there
  is a map
  $$\phi=\phi_{R',k,m',\epsilon'}:X_{R'}(k,m',\epsilon')\to
  Y_R(k,m,\epsilon),$$ which is measure preserving.  We say that $X$
  and $Y$ are asymptotically equal, if both $X$ is asymptotically
  included in $Y$ and $Y$ is asymptotically included in $X$.
\end{definition}

\begin{remark} \label{rmq:asympt}
  Note that if $X$ is asymptotically included into $Y$, we obtain that
  \begin{eqnarray*}
    \sup_R \inf_{m,\epsilon}\limsup_k \alpha_k \log
    \mu^X_{R,k,m,\epsilon}(X_R(k,m,\epsilon)) + a_k && \\ \leq \sup_R
    \inf_{m,\epsilon}\limsup_k
    \alpha_k\log\mu^Y_{R,k,m,\epsilon}(Y_R(k,m,\epsilon)) + a_k, &&
  \end{eqnarray*}
  for all sequences $a_k$, $\alpha_k$.
\end{remark}

It is not hard to see that the sets
$$\Gamma_R(Y_1,\dots,Y_n;k,m,\epsilon)$$ and
$$\Gamma^\sa_R(\Re(Y_1),\Im(Y_1),\dots,
\Re(Y_n),\Im(Y_n);k,m,\epsilon)$$ are asymptotically equal; the
relevant maps $\phi$ send the $n$-tuple $(y_1,\dots,y_n)$ of
non-self-adjoint matrices to the $2n$-tuples of self-adjoint matrices
$(\Re(y_1),\Im(y_1),\dots,\Re(y_n),\Im(y_n))$.  This implies (using
the Remark~\ref{rmq:asympt}) that $$\chi(Y_1,\dots,Y_n) =
\chi^\sa(\Re(Y_1),\Im(Y_1),\dots,\Re(Y_n),\Im(Y_n)).$$

We proceed to prove several lemmas that will be used in the proof of
the main theorem.  

\begin{lemma} \label{lemma:push}
  Let $\Gamma\subset M_k^+$ and $U_k\subset U(k)$ be measurable sets.
  Let $$U_k\Gamma = \{vp: v\in U_k, p\in \Gamma\}\qquad
  \text{and}\qquad S(\Gamma) =  \{ \frac{p^2}{2} : p\in \Gamma\}.$$
  Then  $$\lambda_k(U_k\Gamma) = C_k\,
  \mu_k(U_k)\lambda_k^+(S(\Gamma)).$$ In other words, the map
  $\FF:(v,p)\mapsto v\sqrt{2p}$ from $U(k)\times M_k^+$, endowed
  with the measure $\mu_k\times C_k\lambda_k^+$, to $M_k$, endowed
  with the measure $\lambda_k$, is measure preserving.
\end{lemma}

\begin{proof}
  Since invertible matrices are a set of comeasure zero in $M_k$, we
  see by existence of polar decomposition that $P:(v,p)\mapsto vp$ is
  invertible as a map of measure spaces.  We start by computing the
  pull-back of Lebesgue measure on $M_k$ to $U(k)\times M_k^+$.  Note
  that since $P$ is equivariant with respect to the actions of $U(k)$
  by left multiplication, and Lebesgue measure is invariant under this
  action (since the Euclidean structure is), the resulting measure on
  $U(k)\times M_k^+$ is the product of Haar measure on $U(k)$ and some
  measure $\nu_k$ on $M_k^+$, hence $\lambda_k(U_k \Gamma)=\mu_k(U_k)
  \nu_k(\Gamma)$. It remains to identify $\nu_k$.

  We have the equation
  \begin{equation} \label{eqn:fornu}
    d\mu_k(v) d\nu_k(p)= (P^*(\Omega_k):\Omega^u_k\wedge \Omega^+_k)
    d\mu_k(v) d\lambda^+_k(p),
  \end{equation} 
  where $P^*(\Omega_k):\Omega^u_k\wedge \Omega^+_k$ is the ratio of
  the two volume forms.  Furthermore, in view of the mentioned
  invariance under an action of $U(k)$, it is sufficient to compute
  $(P^*(\Omega_k):\Omega^u_k\wedge \Omega^+_k)$ in (\ref{eqn:fornu})
  at the point $(1,p)\in U(k)\times M_k^+$.

  Note that the tangent space $T_{1,p}(U(k)\times M_k^+)$ is isomorphic
  to the direct sum ${\mathfrak u}(k)\times M_k^\sa$, where ${\mathfrak
    u}(k) = i M_k^\sa$ is the Lie algebra of $U(k)$.  Identify
  $T_{(1,p)}(U(k)\times M_k^+) = iM_k^\sa \oplus M_k^\sa$ with $M_k =
  T_p(M_k)$.  Then the inner product given by the trace $\langle
  a,b\rangle = \Re\Tr(ab^*)$ defines on $T_{1,p}$ a Euclidean structure,
  for which the subspaces $M_k^\sa$ and $iM_k^\sa$ are perpendicular.
  Since the restriction of this inner product to ${\mathfrak u}(k)$ is
  the Killing form on this Lie algebra, and the restriction to $T_p
  M_k^+$ is the inner product we chose before on this space, $\Omega_k$
  (which via the above identification is a volume form on $U(k)\times
  M_k^+$) has the form $C_k \Omega^u_k\wedge \Omega^+_k$.  Further,
  $C_k$ is the ratio of the volume form on $U(k)$ arising from the
  Euclidean structure on ${\mathfrak u}(k)$ coming from the Killing
  form and the volume form corresponding to the normalized Haar measure.
  Hence $C_k$ is just the volume of $U(k)$ with respect to the volume
  form arising from the Euclidean structure on ${\mathfrak u}(k)$ coming
  from the Killing form.

  Thus from (\ref{eqn:fornu}) we get that $$d\nu_k(p)d\mu_k(v) = C_k
  d\mu_k (v) \det(DP)(p) d\lambda_k^+(p).$$ It remains to compute $DP$.
  We note that $P$ is the identity map restricted to $M_k^+$.  Choose a
  basis in which $p$ is diagonal with eigenvalues
  $l_1,\dots,l_k$, and let $e_{ij}\in M_k$ be the matrix all
  of whose entries are zero, except that the $i,j$-th entry is $1$.
  Consider the orthonormal basis $\xi_{\alpha\beta}$ for $iM^\sa_k$,
  given by:
  $$ \xi_{\alpha\beta} =
  \begin{cases} \frac{1}{\sqrt{2}}(e_{\alpha
      \beta} - e_{\beta \alpha}) &\text{if }\alpha < \beta \\ 
    ie_{\alpha\alpha}&\text{if } \alpha=\beta \\ i\frac{1}{\sqrt{2}}
    (e_{\alpha\beta} + e_{\beta\alpha}) &\text{if } \alpha > \beta
  \end{cases} 
  $$ Then $$DP(\xi_{\alpha\beta})p = \xi_{\alpha\beta} p = \frac{1}{2}
  (l_\alpha + l_\beta) \xi_{\alpha\beta} +
  \frac{1}{2}(l_\alpha - l_\beta) \eta_{\alpha\beta},\qquad
  \eta_{\alpha\beta}\in M_k^\sa.$$ It follows that $$\det(DP)(p) =
  \frac{1}{2^{k^2}} \prod_{\alpha,\beta=1}^k (l_\alpha +
  l_\beta).$$ Hence we record the final answer:
  $$d\nu_k(p) = C_k 2^{-k^2} \prod_{\alpha,\beta=1}^k (l_\alpha +
  l_\beta) d\lambda_k^+(p)$$ where $l_i$ are the eigenvalues of
  $p$.

  Consider the map $S:p\mapsto \frac{p^2}{2}$ from $M_k^+$ to itself.
  This map is a.e.~invertible; moreover, its Jacobian $\det(DS)$ at
  $p$ is given by $\det(\frac{1}{2}(1\otimes p + p\otimes 1))$, where
  $1\otimes p$ and $p\otimes 1$ are viewed as elements of $M_k \otimes
  M_k \cong M_{k^2}$ (see e.g.  \cite{dvv:entropy2}).  To compute this
  determinant, let $\zeta_i$, $i=1,\dots,k$ be orthonormal
  eigenvectors of $p$, such that $p\zeta_i=l_i\zeta_i$.  Then
  $\zeta_i\otimes \zeta_j$ is an orthonormal basis for ${\mathbb
    C}^{k^2}$, on which $M_{k^2}=M_k\otimes M_k$ acts naturally.
  Moreover, $\frac{1}{2} (1\otimes p + p\otimes 1)(\zeta_i \otimes
  \zeta_j ) = \frac{1}{2}(l_i + l_j) \zeta_i\otimes
  \zeta_j$.  So the determinant is $2^{-k^2}\prod_{\alpha,\beta=1}^k
  (l_\alpha + l_\beta).$ Hence the push-forward of $\nu_k$
  by $S$ is given by $$d(S_*\nu_k)(p) = C_k 2^{-k^2}
  \prod_{\alpha,\beta=1}^k (l_\alpha+l_\beta)
  d\lambda_k^+(p) \cdot \det(DS)^{-1}(p) = C_k d\lambda_k^+(p).$$
  Thus we have
  $$S_*\nu_k= C_k \lambda_k^+,$$ which is our assertion.
\end{proof}

We have the following standard lemma (see \cite{dvv:entropy2}). 

\begin{lemma} \label{lemma:aincl} 
  Let $p$ be a positive element in $M$.  Then the sequences of sets
  $\GammasaR(p, m, k, \epsilon)$ and $\GammasaR(p,m,k,\epsilon)\cap
  M_k^+$, each taken with the measure $\lambda_k$, are asymptotically
  equal.
\end{lemma}

\begin{lemma} \label{lemma:limCk}
  $\lim_k \frac{1}{k^2}\log(C_k) +\frac{1}{2} \log k
  =\frac{3}{4}+\frac{1}{2}\log 2\pi.$
\end{lemma}

In this exact form this lemma can be found, for example, in
\cite{volumes-liegroups} (the reader is cautioned that the cited paper
uses a slightly different normalization of the Killing form, different
from ours by a factor).

\begin{lemma} \label{lemma:ineq}
  Let $y\in (M,\tau)$ be a (not necessarily self-adjoint) random
  variable.  Then $$\chi(y) \leq \chisa\left(\frac{y^*y}{2}\right) +
  \frac{3}{4} + \frac{1}{2}\log 2\pi.$$
\end{lemma}
\begin{proof}
  Denote by $S: M_k \to M_k^+$ the map $$y\mapsto \frac{y^*y}{2}.$$
  Note that $$S(\Gamma_R(y;m,k,\epsilon)) \subset
  \Gamma^\sa_{R^2}\left(\frac{y^*y}{2}; m/2, k, \epsilon\right),$$ hence the
  former is asymptotically included in the latter.  Note that
  $$\Gamma_R(y;m,k,\epsilon) \subset U(k) \Gamma_R(y;m,k,\epsilon).$$
  We therefore get
  \begin{eqnarray*}
    \lambda_k(\Gamma_R(y;m,k,\epsilon)) &\leq&
    \lambda_k(U(k)\Gamma_R(y;m,k,\epsilon)) \\ &\leq&
    \lambda_k(U(k)\{a^*a:a\in \Gamma_R(y;m,k,\epsilon)\})\\ &\leq& C_k
    \lambda_k(S(\Gamma_R(y;m,k,\epsilon))) \\ &\leq&
    C_k\lambda_k\left(\Gamma^\sa_{R^2}\left(\frac{y^*y}{2};\frac{m}{2},k,
    \epsilon\right)\right).
  \end{eqnarray*} 
  Taking the logarithm and passing to the limits gives the result.
\end{proof}

\begin{lemma} \label{lemma:utimes} 
  Let $u,b\in (M,\tau)$ be such that $u$ is a Haar unitary $*$-free
  from the positive element $b$.  Let $z=ub$.  Given $\delta>0$, there
  exists $k_0$, such that for all $k>k_0$, there is a subset
  $X_k\subset U(k)\times \Gamma_R^\sa(\frac{z^*z}{2};m,k,\epsilon)$,
  $$\frac{1}{k^2}\log\frac{\mu_k\times
    \lambda_k^+(X_k)}{\mu_k\times\lambda_k^+(U_k\times
    \Gamma_R^\sa(\frac{z^*z}{2};m,k,\epsilon))\cap M_k^+}\geq
  -\delta,$$ such the map
  \begin{equation} \label{eqn:mapofhypot}
    Q:(v,p) \mapsto v \sqrt{|2p|}
  \end{equation}
  is an asymptotic inclusion of $X_k$, endowed with the measure
  $\mu_k\times C_k \lambda_k^+$, into $\Gamma_R(z;m,k,\epsilon)$,
  endowed with the measure $\lambda_k$.
\end{lemma}
\begin{proof}  Note that by Lemma~\ref{lemma:push}, the map defined in
  equation~(\ref{eqn:mapofhypot}) is measure preserving.  

  Let $R>0$, $\epsilon>0$ and $\delta>0$ be fixed.  By Corollary~2.12
  of \cite{dvv:improvedrandom}, there exists $k_0$, such that for all
  $k>k_0$, and any $x\in M_k^+$, $\|x\|<R$, there is a subset
  $U_k(x)\subset U(k)$ with $\log \mu_k(U_k(x))>-\delta$, so that
  $U_k(x)\cdot x \in\Gamma_R(wx;m,k,\epsilon)$, where $w$ is a Haar
  unitary $*$-free from $x$ (in other words, ``elements of $U_k(x)$ and
  $x$ are $*$-free to order $m$'').  Let
  $$X_k=\bigcup_{x\in\Gamma_R^\sa(\frac{z^*z}{2};m,k,\epsilon)\cap
    M_k^+} U_k(x)\times \{x\}.$$ Since whenever $x\in
  \Gamma_R\left(\frac{z^*z}{2};m,k,\epsilon\right)\cap M_k^+$,
  $U_k(x)\cdot \sqrt{2x} \subset \Gamma_R(z;m,k,\epsilon)$, $Q(X)$
  lies in $\Gamma_R(z;m,k,\epsilon)$. Moreover, since
  $\mu_k(U_k(x))\geq \exp(-\delta)$ for all $x$, we know that the
  volume of $X_k$ with respect to the measure $\mu_k\times\lambda_k^+$
  is at least $\exp(-\delta)$ times that of
  $\Gamma_R(\frac{z^*z}{2};m,k,\epsilon)$.
\end{proof}
\begin{proof}[Proof of \ref{thrm:micro:1} and \ref{thrm:micro:2} %
  in Theorem~\ref{thrm:microstates}] Assume that $x$, $u$ and $b$ are
  as in the statement of Theorem~\ref{thrm:micro:2} and let $z=ub$;
  note that $z$ is $R$-diagonal.  By Lemma~\ref{lemma:utimes} and
  Lemma~\ref{lemma:limCk}, we have that
  $$\chisa(\frac{z^*z}{2}) + \frac{3}{4} + \frac{1}{2}\log 2\pi \leq
  \chi(z).$$  Since, by Lemma~\ref{lemma:ineq}, we always have the
  other inequality, we obtain
  \begin{equation} \label{formulaforchi}
    \chi(z)=\chisa(\frac{z^*z}{2})+\frac{3}{4}+\frac{1}{2}\log 2\pi 
  \end{equation}
  This can be expressed in terms of the free entropy of the symmetric
  variable $x$ as follows (by using the explicit formula for $\chisa$
  of one variable given by Voiculescu in \cite{dvv:entropy2}):
  \begin{eqnarray*}
    \chi(z) &=& \chisa\left(\frac{z^*z}{2}\right) + \frac{3}{4} +
    \frac{1}{2}\log 2\pi \\ &=& 2\left(\frac{3}{4} + \frac{1}{2}\log
    2\pi\right) + \iint \log |s-t| d\mu_{\frac{z^*z}{2}}(s)
    d\mu_{\frac{z^*z}{2}}(t) \\ &=& 2\left( \frac{3}{4} +
    \frac{1}{2}\log 2\pi\right) + 2 \iint \log|s-t| d\mu_{2^{-1/2}
      x}(s)d\mu_{2^{-1/2} x}(t) \\ &=& 2\chisa(2^{-1/2}x).
  \end{eqnarray*}
  This proves
  \ref{thrm:micro:2}.  

  Combining the above with Lemma~\ref{lemma:ineq} we get
  \ref{thrm:micro:1}:
  \begin{eqnarray*}
    \chi(y) &\leq& \chisa\left(\frac{y^*y}{2}\right) + \frac{3}{4} +
    \frac{1}{2}\log 2\pi \\ 
    &=& \chisa\left(\frac{z^*z}{2}\right) + \frac{3}{4} +
    \frac{1}{2}\log 2\pi \\ 
    &=& \chi(z).
  \end{eqnarray*} 
\end{proof}

\begin{proposition} (A change of variables formula for polar
  decomposition) 
  \label{prop:changevar} Let
  $y_1,\dots,y_n$ be elements of a $W^*$-probability space $(M,\tau)$,
  and let $y_i = v_i (y_i^* y_i)^{1/2}$ be their polar decompositions.
  Assume that $f_i: [0,+\infty) \to [0,+\infty)$ are
  $C^1$-diffeomorphisms, and let $z_i = v_i [2f(y_i^*y_i/2)]^{1/2}$.
  Then
  \begin{equation} \label{eqn:changevar}
    \chi(z_1,\dots,z_n) = \chi(y_1,\dots,y_n) + \sum_{j=1}^n \iint
    \log \left| \frac{f(s)-f(t)}{s-t}\right| d\mu_i(s) d\mu_i(t),
  \end{equation}
  where $\mu_i$ is the distribution of $y_i^*y_i/2$ for  $i=1,\dots,n$.
  The same statement holds for $\chi^\omega$ in the place of $\chi$.
\end{proposition}
\begin{proof}
  If for some $i$ the distribution of $y_i^*y_i$ contains atoms, then
  so does the distribution of $z_i^*z_i$.  Indeed, in this case we
  have
  $$\chi(y_1,\dots,y_n) \leq \sum_j \chi(y_j) =-\infty,$$ since by
  Lemma~\ref{lemma:ineq}, $\chi(y_i)\leq
  \chisa(y_i^*y_i/2)+\text{const}=-\infty$.  Similarly,
  $\chi(z_1,\dots,z_n)=-\infty$, and there is nothing to prove.  Hence
  we may assume that the distributions of $y_i^*y_i$, and thus the
  distributions of $z_i^*z_i$
  are non-atomic for all $i$; in particular, that $v_i$ are unitaries.

  We may also assume that $f_i$ for $i\neq 1$ are the identity
  diffeomorphisms; moreover, by replacing $f_i$ with $f_i^{-1}$, we
  only need to prove that the left-hand side of the statement of
  equation~(\ref{eqn:changevar}) is greater than or equal to the right
  hand side.  We write $f=f_1$.

  Consider the mappings $$\FFF: M_k\ni x\mapsto v [2 f(x^*x/2)]^{1/2}\in
  M_k,$$ where $x=v (x^*x)^{1/2}$ is the polar decomposition of $x$,
  and $$\hat{\FFF}: M_k^n\ni (x_1,\dots,x_n) \mapsto (\FFF(x_1),
  x_2,\dots,x_n) \in M_k^n.$$ Note that the set
  $\hat{\FFF}(\Gamma_R(y_1,\dots,y_n;m,k,\epsilon))$, taken with the
  measure $\lambda_k\times\dots\times\lambda_k$ is asymtotically
  included into the set $\Gamma_R(z_1,\dots,z_n;m,k,\epsilon)$, taken
  with the same measure.  Moreover, the infimum of the Jacobian of
  $\hat{\FFF}$ on the set $\Gamma_R(y_1,\dots,y_n;m,k,\epsilon)$ is not
  less than the infimum of the Jacobian of $\FFF$ on the set
  $\Gamma_R(y_1;m,k,\epsilon)$.  View $\FFF$ as a map from $U(k)\times
  M_k^+$ to itself, using the identification of measure spaces
  $U(k)\times M_k^+ \cong M_k$, $(v,p)\mapsto vp$.  Then $\FFF$ acts
  trivially on the unitary component.  Recall that the measure on
  $M_k^+$, arising from the identification of $M_k$ with $U(k)\times
  M_k^+$, is the push-forward of Lebesgue measure on $M_k^+$ to
  $M_k^+$ by the map $p\mapsto p^2/2$.  Hence the infimum of the
  Jacobian of $\FFF$ is equal to the infimum of the Jacobian of the map
  $p\mapsto [2 f(p^2/2)]^{1/2}$ viewed as a map from $M_k^+$ endowed
  with Lebesgue measure to itself, on the set $\Gamma_R(y_1^*y_1/2; m, k,
  \epsilon)$.  The rest of the computation is exactly as in the proof
  of Proposition~3.1 of \cite{dvv:entropy4}.
\end{proof}

\begin{remark} \label{rmq:relchangevar}
  Let $B\subset M$ be a subalgebra of $M$.  The proof of the
  proposition above also works if we replace $\chi(\cdot)$ with the
  relative entropy $\chi(\cdot | B)$ introduced in
  \cite{shlyakht:relmicro}; we leave the details to the reader.
\end{remark}
  
\begin{proof}[Proof of \ref{thrm:micro:3} of Theorem~\ref{thrm:microstates}]
  Assume that $\chi(y) = \chi(ub) > -\infty$.  Because of part
  \ref{thrm:micro:2}, we conclude that $\chi(b) > -\infty$; in
  particular, the distribution of $b$ is non-atomic (see
  \cite{dvv:entropy2}).  Since $(y^*y)^{1/2}=b$, this implies that in
  the polar decomposition of $y = v (y^*y)^{1/2}$, $v$ is a unitary.

  Arguing as in Lemma~4.2 of \cite{dvv:entropy4}, we may assume that
  there exists a family $f_i$ of $C^1$ diffeomorphisms on
  $[0,+\infty)$, and a continuous function $f:[0,+\infty) \to
  [0,+\infty)$, such that $f(\frac{y^*y}{2})$ is the square of a
  $(0,1)$-semicircular random variable, $\|f_j (y^*y) - f(y^*y) \| \to
  0$ as $j\to \infty$, $W^*(y^*y) = W^*(f(y^*y))$, and $\lim_j
  \chi^\sa(f_j(y^*y)) = \chi^\sa (f(y^*y))$.  Let $y=v(y^*y)^{1/2}$ be
  the polar decomposition of $y$; let $z = v [2f({y^*y}/2 )]^{1/2}$,
  and similarly $z_j = v[2 f_j(y^*y/2)]^{1/2}$.  Then by
  Proposition~\ref{prop:changevar} and the explicit formula for the
  free entropy of one variable given by Voiculescu (Proposition~4.5 in
  \cite{dvv:entropy2}), we get for all $j$,
  $$\chi(z_j) = \chi(y) + \chi^\sa\left(f_j\left(\frac{y^*y}{2}\right)\right) -
  \chi^\sa\left(\frac{y^*y}{2}\right).$$ Applying Proposition~2.6 of
  \cite{dvv:entropy2}, we get that
  \begin{eqnarray*}
    \chi(z) &\geq & \limsup_j \chi(z_j)\\ & = & \limsup_j \Bigg[
    \chi(y) + \chi^\sa\left(f_j\left(\frac{y^*y}{2}\right)\right) -
    \chi^\sa\left(\frac{y^*y}{2}\right)\Bigg] \\ &=& \chi(y) +
    \chi^\sa\left(f\left(\frac{y^*y}{2}\right)\right) -
    \chi^\sa\left(\frac{y^*y}{2}\right).
  \end{eqnarray*}
  Since $\chi(y) = \chi(ub)$ by assumption, and $\chi(ub) = \chi^\sa(
  \frac{y^*y}{2}) + \frac{3}{4} + \frac{1}{2}\log 2\pi$ by
  Theorem~\ref{thrm:micro:2} we get that $$\chi(z) \geq
  \chi^\sa\left(f\left(\frac{y^*y}{2}\right)\right) + \frac{3}{4} +
  \frac{1}{2}\log 2\pi.$$ By assumption, the distribution of
  $(z^*z)^{1/2}$ is quarter-circular (i.e., it is the absolute value
  of a $(0,2)$-semicircular).  Let $c$ be a circular variable (i.e.,
  its real and imaginary parts are free $(0,1)$-semicircular
  variables).  Then, since $c$ is $R$-diagonal (see
  \cite{speicher:r-diag}), we have by \ref{thrm:micro:2}, that
  \begin{eqnarray*}
    \chi(c) &=& \chi^\sa\left( \frac{c^*c}{2}\right) + \frac{3}{4} +
    \frac{1}{2}\log 2\pi\\&=& \chi^\sa\left(
    f\left(\frac{y^*y}{2}\right)\right) + \frac{3}{4} +
    \frac{1}{2}\log 2\pi,
  \end{eqnarray*}
  since $c^*c$ has the same distribution as $z^*z=2f({y^*y}/2)$.  
  Hence $\chi(z)\geq \chi(c)$.  

  On the other hand, $c$ is $R$-diagonal, with the same distribution
  of the positive part as $z$, so by \ref{thrm:micro:1}, we have
  $\chi(z) \leq \chi(c)$.  So $\chi(z)= \chi(c)$.

  We claim that $z$ is circular.  This will prove the proposition,
  since then the polar and positive parts of $z$ are $*$-free (see
  \cite{dvv:circular} or \cite{speicher:r-diag}), and thus the polar
  and positive parts of $y$ are $*$-free, since the polar part of $y$
  is the same as the polar part of $z$, and the positive part of
  $y$ is some function of the positive part of $z$.

  Now, for the claim that $z$ is circular, let $\gamma$ be a
  complex number of modulus one; then $\chi(\gamma z) =
  \chi(z)$.  Let $$X_\gamma = \frac{1}{2}( \gamma z +
  \overline{\gamma} z^*), \quad Y_\gamma=\frac{1}{2i}(\gamma z -
  \overline{\gamma} z^*).$$ Then $$\tau(X_\gamma^2) =
  \frac{1}{4}\left[ 2 \tau(z z^*) + \gamma^2 \tau(z^2) +
  \overline{\gamma^2}\cdot \overline{\tau(z^2)} \right] {} .$$ Similarly,
  $$\tau(Y_\gamma^2) = \frac{1}{4} \left[ 2 \tau(zz^*) -
  \gamma^2\tau(z^2) - \overline{\gamma^2} \cdot
  \overline{\tau(z^2)}\right]{}.$$ We choose $\gamma$ such that
  $\gamma^2 \tau(z^2)$ is purely imaginary.  Since $\tau(z^*z)=2$, we
  have then $\tau(X_\gamma^2) = \tau(Y_\gamma^2)=1$.  But $\chi(z) =
  \chi(c)=\chi^\sa(x_1, x_2)$, where $x_i$ are free $(0,1)$
  semicircular variables.  Hence we have 
  $$\chi(z) = \chi^\sa (X_\gamma, Y_\gamma) = \chi(\gamma z) =
  \chi^\sa(x_1,x_2),$$ where $X_\gamma$ and $Y_\gamma$ are some
  self-adjoint random variables of covariance $1$.  But then by
  Voiculescu's Proposition~2.4 of \cite{dvv:entropy4}, $X_\gamma$ and
  $Y_\gamma$ are both semicircular and free, so that $\gamma z$ is
  circular, so $z$ is circular.
\end{proof}    

\section{Maximization of free entropy for matrices.}
\begin{theorem}
  Let $X_{ij}$, $1\leq i,j\leq d$ be non-commutative random variables
  in a tracial non-commutative probability space $(M,\hat\tau)$.  Let
  $Z\in M\otimes M_d$ be given by $$Z=\sum_{i,j=1}^d X_{ij}\otimes
  e_{ij},$$ where $e_{ij}$ are matrix units in the algebra of $d\times
  d$ matrices. We denote by $\tau$ the normalized trace on $M\otimes
  M_d$.  Let $\omega$ be a free ultrafilter.  Let $X$ be a
  self-adjoint variable with $\tau(X^{2n+1})=0$ for all $n\in \mathbb
  N$, and such that $\tau(X^{2n}) =\tau((Z^*Z)^n)$, $\forall
  n\in\mathbb N$.  Then we have
  \begin{enumerate}
    \renewcommand{\theenumi}%
    {\arabic{section}.\arabic{subsection}(\alph{enumi})}
    \renewcommand{\labelenumi}{(\alph{enumi})}
  \item $\chi^{d\omega}(\{X_{ij}\}_{1\leq i,j\leq d})\leq d^2
    \chi^\omega(Z) + d^2\log d\leq 2d^2 \chi^\sa(2^{-\frac{1}{2}}X) +
    d^2 \log d$.
    \label{thrm:matrix:1}
  \item Equality holds in \ref{thrm:matrix:1} if $Z$ is $R$-diagonal
    and $*$-free from the algebra $1\otimes M_d$. \label{thrm:matrix:2}
  \item If equality holds in \ref{thrm:matrix:1} and
    $\chi^\sa(2^{-\frac{1}{2}}X)\neq -\infty$, then $Z$ is $R$-diagonal
    and is $*$-free from the algebra $1\otimes M_d$. \label{thrm:matrix:3}
  \end{enumerate}
\end{theorem}
\begin{proof}
  Let $B=1\otimes M_d$.  We have by \cite{shlyakht:relmicro} that
  $$\chi^{d\omega}(\{X_{ij}\}) = d^2 \chi^\omega(Z|B)+d^2\log d \leq
  d^2 \chi^\omega(Z) + d^2\log d.$$ (We have the summand $d^2\log d$
  rather than $\frac{d^2}{2}\log d$ appearing above because we are
  dealing with $\chi$, not $\chi^\sa$).  Moreover, $d^2 \chi^\omega(Z)
  \leq 2d^2\chi^\sa(2^{-\frac{1}{2}}X)$ by
  Theorem~\ref{thrm:microstates}, hence \ref{thrm:matrix:1}.
  
  If $Z$ is $*$-free from $B$, then, by \cite{shlyakht:relmicro}, we
  have $\chi^\omega(Z|B)=\chi^\omega(Z)$.  Moreover, if $Z$ is
  $R$-diagonal, we have, by Theorem~\ref{thrm:microstates}, that
  $\chi^\omega(Z) =2\chi^\sa(X/\sqrt 2)$, which proves
  \ref{thrm:matrix:2}.

  Assuming the conditions in \ref{thrm:matrix:3} are satisfied, we get
  that $\chi^\omega(Z)= 2\chi^\sa(2^{-\frac{1}{2}}X) >-\infty$, so $Z$
  is $R$-diagonal by Theorem~\ref{thrm:micro:3}, i.e, $Z$ has polar
  decomposition $Z=v(Z^*Z)^{1/2}$, where $v$ is a Haar unitary, which
  is $*$-free from $Z^*Z$. Note also that we are given that
  $\chi^\omega(Z|B)=\chi^\omega(Z)$. We may assume, as in the proof of
  statement \ref{thrm:micro:3} of Theorem~\ref{thrm:microstates} that
  there exists a family $f_i$ of $C^1$ diffeomorphisms on
  $[0,+\infty)$, and a continuous function $f:[0,+\infty)\to
  [0,+\infty)$, such that $f(\frac{Z^*Z}{2})$ is the square of a
  $(0,1)$-semicircular random variable, $\|f_j (Z^*Z) - f(Z^*Z) \| \to
  0$ as $j\to \infty$, $W^*(Z^*Z) = W^*(f(Z^*Z))$, and $\lim_j
  \chi^\sa(f_j(Z^*Z)) = \chi^\sa (f(Z^*Z))$.  Given the polar
  decomposition $Z=v(Z^*Z)^{1/2}$, let $z = v [2f({Z^*Z}/2 )]^{1/2}$,
  and similarly $z_j = v[2 f_j(Z^*Z/2)]^{1/2}$.  Notice that $z$ is
  circular; moreover, since $W^*(Z^*Z)=W^*(f(Z^*Z))=W^*(z^*z)$, we
  have that $Z\in W^*(z)$.  Hence it will suffice to prove that $z$ is
  $*$-free from $B$, as then also $Z$ is $*$-free from $B$.

  By Remark~\ref{rmq:relchangevar} and the
  explicit formula for the free entropy of one variable given by
  Voiculescu (Proposition~4.5 in \cite{dvv:entropy2}), we get for all
  $j$,
  $$\chi^\omega(z_j|B) = \chi^\omega(Z|B) +
  \chi^\sa\left(f_j\Big(\frac{Z^*Z}{2}\Big)\right) -
  \chi^\sa\left(\frac{Z^*Z}{2}\right).$$ We get
  \begin{eqnarray*}
    \chi^\omega(z|B) &\geq & \limsup_j \chi^\omega(z_j|B)\\ & = &
    \limsup_j \Bigg[ \chi^\omega(Z|B) +
    \chi^\sa\left(f_j\Big(\frac{Z^*Z}{2}\Big)\right) -
    \chi^\sa\left(\frac{Z^*Z}{2}\right)\Bigg] \\ &=& \chi^\omega(Z|B) +
    \chi^\sa\left(f\Big(\frac{Z^*Z}{2}\Big)\right) -
    \chi^\sa\left(\frac{Z^*Z}{2}\right).
  \end{eqnarray*}
  By assumption, we have that $\chi^\omega(Z|B)=\chi^\omega(Z)$;
  moreover, by $R$-diagonality of $Z$ we get by Theorem~\ref{thrm:micro:2} that
  $\chi(Z)=\chisa(Z^*Z/2) + 3/4 + (1/2)\log 2\pi$.
  Therefore, we get that
  \begin{eqnarray*}
    \chi^\omega(z|B) &\geq& \chisa\left(\frac{Z^*Z}{2}\right) +
    \frac{3}{4} + \frac{1}{2}\log 2\pi \\ &&
    +\chi^\sa\left(f\Big(\frac{Z^*Z}{2}\Big)\right) -
    \chi^\sa\left(\frac{Z^*Z}{2}\right) \\ &=&\frac{3}{4} +
    \frac{1}{2}\log
    2\pi+\chi^\sa\left(f\Big(\frac{Z^*Z}{2}\Big)\right).
  \end{eqnarray*}
  But $z$ is circular, in particular $R$-diagonal; moreover,
  $z^*z/2=f(Z^*Z/2)$.  So from the formula in
  \ref{thrm:micro:2}, we get that 
  $$\chi^\omega(z)=\frac{3}{4} + \frac{1}{2}\log
  2\pi+\chi^\sa\left(f\Big(\frac{Z^*Z}{2}\Big)\right).$$ Thus
  $\chi^\omega(z|B)\geq \chi^\omega(z)$.  Since $\chi^\omega(z|B)\leq
  \chi^\omega(z)$ in general, we get that
  $\chi^\omega(z|B)=\chi^\omega(z)$.

  Now let $S_1$, $S_2$ be the real and imaginary parts of $z$.  Then
  we have that $\chi^\sa(S_1,S_2|B)=\chi^\sa(S_1,S_2)$.  Since $S_1$
  and $S_2$ are two free semicircular variables, it follows by
  Theorem~4.5 from \cite{shlyakht:relmicro} that $W^*(S_1,S_2)$ is
  free from $B$.  Hence $z$ is $*$-free from $B$; hence $Z$ is
  $*$-free from $B$.
\end{proof}
\bibliographystyle{amsplain}

\begin{thebibliography}{10}

\bibitem{hiai-petz}
F.~Hiai and D.~Petz, \emph{Properties of free entropy related to polar
  decomposition}, Preprint, 1998.

\bibitem{volumes-liegroups}
M.~S. Marinov, \emph{Invariant volumes of compact groups}, J. Phys. A
  \textbf{13} (1980), 3357--3366.

\bibitem{nss:entropy}
A.~Nica, D.~Shlyakhtenko, and R.~Speicher, \emph{Some minimization problems for
  the free analogue of the fisher information}, Preprint, 1998.

\bibitem{speicher:r-diag}
A.~Nica and R.~Speicher, \emph{{$R$}-diagonal pairs---a common approach to
  {Haar} unitaries and circular elements}, Free Probability (D.-V. Voiculescu,
  ed.), 1997, pp.~149--188.

\bibitem{shlyakht:relmicro}
D.~Shlyakhtenko, \emph{A microstates approach to relative free entropy},
  Preprint, 1998.

\bibitem{shlyakht:cpentropy}
\bysame, \emph{Free entropy with respect to a completely-positive map},
  Preprint, 1998.

\bibitem{dvv:circular}
D.-V. Voiculescu, \emph{Circular and semicircular systems and free product
  factors}, Operator Algebras, Unitary Representations, Enveloping Algebras,
  and Invariant Theory, Progress in Mathematics, vol.~92, Birkh\"auser, Boston,
  1990, pp.~45--60.

\bibitem{dvv:entropy2}
\bysame, \emph{The analogues of entropy and of {Fisher's} information measure
  in free probability theory {II}}, Invent. Math. \textbf{118} (1994),
  411--440.

\bibitem{dvv:entropy4}
\bysame, \emph{The analogues of entropy and of {Fisher}'s information measure
  in free probability theory, {IV}: Maximum entropy and freeness}, Free
  Probability (D.-V. Voiculescu, ed.), American Mathematical Society, 1997,
  pp.~293--302.

\bibitem{dvv:entropy5}
\bysame, \emph{The analogues of entropy and of {Fisher}'s information measure
  in free probabilility, {V}}, Invent. Math. \textbf{132} (1998), 189--227.

\bibitem{dvv:improvedrandom}
\bysame, \emph{A strengthened asymptotic freeness result for random matrices
  with applications to free entropy}, IMRN \textbf{1} (1998), 41 -- 64.

\end{thebibliography}

\providecommand{\bysame}{\leavevmode\hbox to3em{\hrulefill}\thinspace}

\end{document}